# Online Allocation Rules in Display Advertising

Davood Shamsi

Department of Management Science and Engineering, Stanford University, CA 94305,
AOL Platforms, Palo Alto, CA, 94306, davood@stanford.edu

Marius Holtan

AOL Platforms, Palo Alto, CA, 94306, marius.holtan@teamaol.com

Robert Luenberger

AOL Platforms, Palo Alto, CA, 94306, rob.luenberger@teamaol.com

Yinyu Ye

Department of Management Science and Engineering, Stanford University, CA 94305, yinyu-ye@stanford.edu

Efficient allocation of impressions to advertisers in display advertising has a significant impact on advertisers' utility and the browsing experience of users. The problem becomes particularly challenging in the presence of advertisers with limited budgets as this creates a complex interaction among advertisers in the optimal impression assignment. In this paper, we study online impression allocation in display advertising with budgeted advertisers. That is, upon arrival of each impression, cost and revenue vectors are revealed and the impression should be assigned to an advertiser almost immediately. Without any assumption on the distribution/arrival of impressions, we propose a framework to capture the risk to the ad network for each possible allocation; impressions are allocated to advertisers such that the risk of ad network is minimized. In practice, this translates to starting with an initial estimate of dual prices and updating them according to the belief of the ad network toward the future demand and remaining budgets. We apply our algorithms to a real data set, and we empirically show that Kullback-Leibler divergence risk measure has the best performance in terms of revenue and balanced budget delivery.

*Key words*: online algorithms, linear programming, primal-dual, dynamic price update, risk minimization.

*History*: This paper was first submitted on July 21, 2014.

## 1. Introduction

Online Display Advertising, which includes text ads and banner ads, is a large, growing, dynamic and complex market. US advertisers spent more than $17 billion on display advertising and, as they realize the effectiveness of display ads, the advertisers are spending more advertising dollars on

them. It is expected that display ad spending will surpass $21 billion in 2014. This would represent 44% of all digital advertising, second only to search advertising. In 2015, display advertising may even outstrip search advertising (eMarketer.com 2014). Another aspect of the display advertising landscape is complexity. This includes attribution models (Shao and Li 2011), competition among ad networks for attribution (Berman 2013), sequential auctions (Weber 1981), and advertiser with limited budgets (Goel et al. 2009).

While publishers such as the huffingtonpost.com provide a medium for advertising, ad networks optimize ad placement on these digital ad mediums. For this purpose, ad networks track user activities online and historical interactions between users and ads, typically via tracking pixels. Ad networks use this information to serve users relevant advertisements and help advertisers send their messages to the desired users.

When a user browses a webpage with an advertising slot, an impression request is sent to the ad network. Upon arrival of the impression request, the ad network must assign it to an advertiser (media) almost immediately. In the allocation process, the ad network retrieves historical information about the user and calculates value (bid) of serving each media (based on the estimated probability of a click or other action, e.g.). In the presence of budget constraints for advertisers, allocating the impression to the advertiser with the highest value might be suboptimal. Considering the budget of each advertiser, the impression request should be allocated to an advertiser such that the long term revenue of the ad network is maximized.

The display advertising problem has many challenging aspects. First, budget constraints introduce interdependency among advertisers. Even the offline version of the problem is NP-hard (Goel et al. 2009). However, when the ratio of bids to budgets is small enough, an offline LP formulation of the problem results in a $(1-\epsilon)$-optimal solution (Bahl et al. 2007). In the last few years, many dual-based heuristics have been developed to find near optimal solutions for the online version of this problem (for example Agrawal et al. (2009), Devanur et al. (2011)). Second, even if all impression requests are available in advance, solving the LP to find a $(1-\epsilon)$-optimal solution is not computationally practical. Large ad networks serve billions of impressions per day and currently there is no LP solver than can handle such a large problem in the time allotted. The allocation algorithm should be very fast, highly scalable and robust to outliers.

The third challenge is the unpredictable and unstructured nature of impression arrivals. The stream of impression requests is not stationary and doesn't follow any predefined pattern. For example, an unpredictable breaking news story could flood the ad network with impression requests from news websites. Finally, because of the technological challenges as well as the attribution model being used, the value of a particular allocation (and associated budget consumption) can be observed only with some delay.

In the display advertising allocation problem, there are $N$ advertisers each with budget $b_i$, $i = 1, 2, \ldots, N$, and there are $M$ impressions that should be served. Impressions arrive sequentially and upon arrival they must be assigned to an advertiser. The value of assigning impression $j$ to advertiser $i$ is $r_{ji}$ and this would consume $a_{ji}$ of the $i$-th advertiser's budget. The objective is to maximize the total value of the allocation without violating any budget constraint.

If one had all the impressions beforehand, the display advertising problem could be formulated as a mixed integer linear programming. Let $\mathbf{x}_j = (x_{j1}, x_{j2}, \ldots, x_{jN})$ be the variable that represents allocation of the $j$-th impression, i.e., $x_{ji} = 1$ if $j$-th impression is assigned to the $i$-th advertiser and $x_{ji} = 0$ otherwise. We can find the optimal impression allocation by solving the following maximization problem:

$$\max \sum_{j=1}^{M} \sum_{i=1}^{N} r_{ji} x_{ji} \tag{1}$$

$$\text{s.t.} \sum_{j=1}^{M} a_{ji} x_{ji} \leq b_i \ \forall i$$

$$\sum_{i=1}^{N} x_{ji} \leq 1 \ \forall j, x_{ji} \in \{0, 1\} \ \forall (j, i).$$

Here, $M$ is total number of impressions. In the rest of paper, we denote vectors of revenues and costs by $\mathbf{r}_j = (r_{j1}, r_{j2}, \ldots, r_{jN})$ and $\mathbf{a}_j = (a_{j1}, a_{j2}, \ldots, a_{jN})$, respectively. We can relax the constraint $x_{ji} \in \{0, 1\}$ to $x_{ji} \geq 0$ and turn this into a linear programming problem.

The dual of the linear program would be

$$\min \ \sum_{j=1}^{M} \hat{p}_j + \sum_{i=1}^{N} b_i p_i \tag{2}$$

$$\text{s.t.} \ \hat{p}_j + a_{ji} p_i \geq r_{ji} \ \leftarrow x_{ji} \ \forall (j,i),$$

$$p_i \geq 0 \ \forall i,$$

$$\hat{p}_j \geq 0 \ \forall j.$$

One simple observation is that, for every $j$ at optimality,

$$\hat{p}_j = \max_i \{r_{ji} - a_{ji} p_i\}.$$

Thus, once optimal dual $p_i$'s are found, optimal $\hat{p}_j$ can be constructed sequentially.

A practical, viable solution to the online display advertising problem should allocate the arriving impressions to advertisers almost immediately. When a user visits a particular website containing an advertisement slot, an impression request is sent to the Adserver. The Adserver should respond in less than 5 milliseconds. Thus, a computationally complex algorithm that can't meet 5 milliseconds constraint is not practical.

The online resource allocation problem of (1) has been studied under several assumptions on revenue and cost vectors $(\mathbf{r}_j, \mathbf{a}_j)$. These assumptions can be divided into four main categories.

- Adversarial model: At each time step, an adversarial agent chooses $\mathbf{r}_j$ and $\mathbf{a}_j$ (Mehta et al. 2005, Aggarwal et al. 2011, Devanur et al. 2013, Karp et al. 1990).

- Random permutation model: $\mathbf{r}_j$ and $\mathbf{a}_j$ are fixed at the beginning and arrive according to a random permutation (Agrawal et al. 2009, Devanur and Hayes 2009, Mahdian 2011, Karande et al. 2011, Goel and Mehta 2008, Bhalgat et al. 2012).

- Unknown iid distribution: $\mathbf{r}_j$ and $\mathbf{a}_j$ are independently sampled from an unknown distribution (Devanur et al. 2011, Balseiro et al. 2012, Ghosh et al. 2009, Devanur et al. 2012).

- Known iid distribution: $\mathbf{r}_j$ and $\mathbf{a}_j$ are independently sampled from a known distribution (Feldman et al. 2009b, Jaillet and Lu 2012, Alaei et al. 2012, Manshadi et al. 2012).

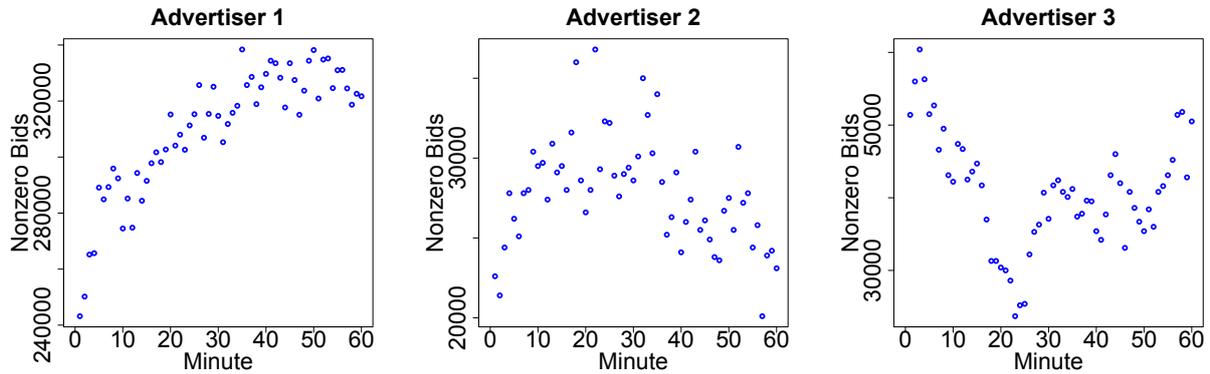

**Figure 1** Number of nonzero bids for three advertisers during an hour. Bids and user types can change by time.

The adversarial model is too restrictive and usually imposes unnecessary constraints on the problem. Since an adversarial model is a worst cases analysis, one does not expect a near optimal algorithm (Mehta et al. 2005). On the other hand, the last three models: random permutation, unknown iid distribution and especially known iid distribution, may assume too much structure about the problem. In practice, the cost and revenue vectors don't arrive in adversarial nor uniform manner. Figure (1) shows number of nonzero bids of three advertisers in one hour (9am-10am PST). The advertiser associated with the left plot (advertiser 1) would have more opportunities as we get closer to the end of the hour. This might be due to the time of day pattern of the Internet users. However, in the same hour, advertisers associated with the middle and right plot of Figure (1) have different opportunity (demand) patterns.

Any reasonable estimation of dual prices associated to Problem (1) can be useful in guiding online allocation. Consider the following small version of Problem (1) using only $\epsilon$-fraction of impressions.

$$\max \ \sum_{j=1}^{\epsilon \cdot M} \sum_{i=1}^{N} r_{ji} x_{ji} \tag{3}$$

$$\text{s.t.} \ \sum_{j=1}^{\epsilon \cdot M} a_{ji} x_{ji} \leq b_i^{\epsilon}, i = 1, 2, \ldots, N,$$

$$\sum_{i=1}^{N} x_{ji} \leq 1, j = 1, 2, \ldots, \epsilon \cdot M, x_{ji} \in \{0, 1\} \ \forall (j, i).$$

Note that we used $\epsilon$-fraction of impressions, and thus budgets should be scaled down accordingly to some $b_i^\epsilon \leq \epsilon b_i$ for all $i$. The dual of the sample linear program would be

$$\min \sum_{j=1}^{\epsilon \cdot M} \hat{p}_j + \sum_{i=1}^{N} b_i^\epsilon p_i \qquad (4)$$
$$\text{s.t.} \quad \hat{p}_j + a_{ji} p_i \geq r_{ji} \leftarrow x_{ji} \; \forall (j,i),$$
$$p_i \geq 0, i = 1, 2, \ldots, N,$$
$$\hat{p}_j \geq 0, j = 1, 2, \ldots, \epsilon \cdot M.$$

Let $p_i^\epsilon$, $i=1,...,N$, be optimal dual prices of the sample linear program. Then, the online allocation for the rest of the impressions, say $j$, becomes: if $\max_i \{r_{ji} - a_{ji} p_i^\epsilon\} \leq 0$ then $x_{ji} = 0$ for all $i$; otherwise select one $i^* = \arg\max\{\max_i \{r_{ji} - a_{ji} p_i^\epsilon\}\}$ and set $x_{ji^*} = 1$.

Let $\gamma$ denote the maximum ratio of budget consumed by impressions to the total budget, i.e.,

$$\gamma = \max_{j,i} \frac{a_{ji}}{b_i}.$$

Under the random permutation model, Devanur and Hayes (2009) and Agrawal et al. (2009) show if set $b_i^\epsilon = \epsilon(1-\epsilon) b_i$ for all $i$ and $\gamma \leq O\left(\frac{\epsilon^3}{N \log(M/\epsilon)}\right)$, such an online allocation algorithm, on expectation, achieves $1 - O(\epsilon)$ fraction of the offline optimal revenue. Moreover, Agrawal et al. (2009) show if the sample optimal dual prices are updated in subsequent doubling intervals, the condition on $\gamma$ can be loosened to $\gamma \leq O\left(\frac{\epsilon^2}{N \log(M/\epsilon)}\right)$.

As mentioned earlier, the random permutation model may not be realistic in practice. Furthermore, the dual prices used in these online allocations are not updated later to reflect the remaining budgets. In this paper, we use the sample linear program optimal dual prices as *initial prices* and propose simple price updating rules that incorporate belief of the ad network regarding the future demand and remaining budgets. To appropriately model the uncertainty of future demand, we formulate the resource allocation problem (1) as a risk minimization problem. At each step, there are uncertainties about the future impression requests, thus, keeping or consuming any resource (budget) would expose us to some risk. In other words, if we keep a resource and do not consume it now, demand for that resource might be very low in the future. Conversely, if we allocate an impression to an advertiser and consume its budget, there might be a high demand for its budget in the future and we lose the opportunity to serve impressions later with higher gain. Thus, upon arrival of each impression request, it is assigned to an advertiser such that the risk of both under and over budget consumption is minimized. Our contributions can be summarized as follows.

**Resource Allocation as a Risk Minimization Problem:** Without assuming any particular structure on arrival of impression requests, we pose the underlying resource allocation problem as a convex risk minimization model problem. The ad network has an initial belief about the future demand and this initial belief is updated as new requests are observed. Using a robust representation of convex risk measures, we show that each distance function in a measurable space of possible outcomes translates to an allocation algorithm.

**Risk Minimization Algorithm:** We show that for each arrival impression, the convex risk minimization model can be solved in a close form so that the allocation can be decided very efficiently. Thus, our risk minimization model is tractable for practical real time allocations.

**Performance on a Real Dataset:** We use real advertisers and impression requests from AOL/Advertising.com to validate the effective performance of our model/algorithm. Among the proposed distance functions, we show that Kullback- Leibler divergence risk measure results in the highest revenue.

**Frequency Capping:** Typically, advertisers limit the number of times a particular ad can be served to a particular user over some time frame. For example, they require the ad network to show at most one ad per hour to a user. Incorporating these constraints would significantly increase the computational complexity of the LP in the initial phase. We introduce a relaxation to the original frequency capped LP that extremely reduces the size of the problem.

Before presenting the outline of the paper, let's remove the ambiguity in four terms that are used throughout the paper. We use the following group of words interchangeably: (1) Impression request, impression, demand, (2) Advertiser, campaign, media, (3) Bid, impression value, (4) Budget, resource.

The rest of the paper is organized as follows. Related work and literature review is presented in Section 2. Section 3 details the risk measure minimization framework for resource allocation. In Section 4, we pose the display advertising problem as a convex risk minimization problem. Then, we study robust representation of the risk measure and build a bridge between the space of value functions over budgets and the space of penalty functions over probability measures on possible future outcomes. In Section 5, we show allocation based on risk minimization has a simple closed form solution and can be quickly computed by ad servers. Finally, results of the allocation algorithms on real data are presented in Section 6.

## 2. Related Work

In an original paper, Mehta et al. (2005) introduced the adwords problem and proposed two allocation algorithms that are based on RANKING and BALANCED algorithms for the online matching problem. The two algorithms have competitive ratios of $1 - 1/e$ which is the theoretical performance bound for the adversarial model. Thus, the algorithms are optimal under the adversarial model. Later, motivated by the display advertising allocation, Feldman et al. (2009a) studied the matching problem under the free disposal assumption. In their model, it is possible to assign more than the desired impressions to the advertisers, and excessive impressions are ignored. They introduced a primal-dual algorithm with exponential updates. The proposed algorithm is $(1 - 1/e)$-competitive which is again theoretically tight under the adversarial model.

Our work is closely related to the work by Devanur and Hayes (2009) on online algorithms for the adwords problem and Agrawal et al. (2009) for general resource allocation problems, both under the random permutation model. We described their results in the previous section, which uses the optimal dual price of each budget constraint of a sampled linear program. In particular, Agrawal et al. (2009) proved that if the dual prices are updated in subsequent doubling intervals, the condition to achieve $(1-\epsilon)$-competitive ratio can be relaxed to $\gamma \leq O\left(\frac{\epsilon^2}{N \log(M/\epsilon)}\right)$, where $\gamma$ is the bid to budget ratio defined early. Furthermore, by constructing an example, they showed that any $(1-\epsilon)$-competitive online algorithm needs $\gamma \leq O\left(\frac{\epsilon^2}{\log N}\right)$. One can see there is a gap between the sufficient and necessary conditions on $\gamma$.

More recently, using the primal solution, Kesselheim et al. (2013) closed the gap. Upon arrival of a demand (impression) $j$, the primal sample linear program, including all early demand information, of the original one is solved. Then, the allocation of the newly arrived impression is decided by a randomized rounding of the primal solution value $x_{ji}$. This algorithm achieves a $(1-\epsilon)$-competitive solution when $\gamma$ is upper bounded by $O\left(\frac{\epsilon^2}{\log d}\right)$; where $d$ is the column sparsity of the demands and equals to $N$ in the worst case. In display advertising with millions of impressions, solving a linear program in each stage can soon becomes a bottleneck. Thus, even though their algorithm touches the optimal bound for $\gamma$, it may not be practical in display advertising, which requires an almost spontaneous allocation.

Also in another closely related work, Devanur et al. (2011) turned to the unknown i.i.d model for inputs, and improved the initial performance bound. They developed heuristic allocation algorithms for online mixed covering-packing and general online resource allocation problems. Taking advantage of the stronger assumed mathematical model, they showed that their resource allocation algorithm is $(1-\epsilon)$-competitive if $\gamma$ is at most $O\left(\frac{\epsilon^2}{\log(N/\epsilon)}\right)$. Consequently, limiting the focus to the adwords problem, Devanur et al. (2012) introduced an asymptotically optimal allocation algorithm under the unknown i.i.d. model.

Ciocan and Farias (2012) proposed a dual-based algorithm for the display advertising allocation under high dimensional demand model. They assumed each impression has $d$ possible features (such as geo, gender,...) that characterize the value of the impression for advertisers. Each advertiser is only interested in a subset of impressions with some specific properties, for example, an advertiser might only be interested in showing ads to male users in California. Advertisers have the same value for all impressions in their target group. Since targets of advertisers can overlap, when the demand dimension $d$ is large, predicting demand (impression) for each advertiser would be intractable. Assuming impressions arrive i.i.d. from an unknown distribution, they showed their algorithm is $(1-\epsilon)$-optimal for small enough bids.

Bhalgat et al. (2012) focused on the smoothness of delivery, which is critical in display advertising. Specifically, they partitioned the flight of the algorithm into many intervals, and then, appropriate constraints were added to the problem to enforce smooth delivery. Using the LP formulation, they developed a primal-dual algorithm and proved performance of the algorithm is tight under the adversarial model. What uniquely distinguishes this work from others is applying the algorithm to a real data set from DoubleClick. They tested the proposed algorithms on 10 data sets and they observed that the performance of an algorithm can be significantly different across different data sets.

As described earlier, upon arrival of each impression request, the online allocation algorithm developed in the current paper quickly solves a risk minimization problem for the allocation decision. Thus, we essentially update the prices of budgets on each arrival. From this perspective, our work is also closely related to the research line of the prediction market design for contingent claims; see Agrawal et al. (2011). Even though prediction markets are designed to aggregate information of the traders, from the market maker's point of view, it is basically a risk minimization problem. The market maker has an initial belief about the event of interest, and claims arrive sequentially. Upon arrival of the claim, the market maker allocates the claim such that the exposed future risk is minimized. And consequently, the market maker's belief (price) is updated. One can think of online resource allocation as an information market, and the market maker (ad network) has an initial belief (price) about the future demand. At each step, a demand vector (impression) arrives and it should be assigned such that the risk is minimized.

The display advertising problem and the adwords problem are studied in other contexts which are not the focus of this work. For example, this problem can be studied as a second price auction which is particularly important with the rise of ad exchanges (like DoubleClick). Sequence of second price auctions is not truthful and strategic behavior might arise (Charles et al. 2013, Leme et al. 2012, Iyer et al. 2011). The display advertising problem is also studied under the assumption of uncertain information (Mahdian et al. 2012) and unknown/volatile demand (Ciocan and Farias 2012).

# 3. Risk Minimization for Display Advertising

In this section, we first define the concept of risk measure in resource allocation. This is, we define a function that assigns a monetary risk value to each state in the process of online resource allocation. Later, in Section 3.2, this concept becomes the foundation for the allocation and price updating rules. For more details on risk measure please see Föllmera and Schied (2002).

## 3.1. Risk Measure Framework for Display Advertising

In order to formally define the risk of a random return, let measurable space $(\Omega, \mathcal{F})$ be the set of possible scenarios, and $\mathcal{X}$ be the set of all bounded measurable functions on $(\Omega, \mathcal{F})$.

$$\mathcal{X} = \{X : \Omega \to \mathbb{R} \mid X \text{ is bounded measurable functions on}(\Omega, \mathcal{F})\}.$$

Functional $X \in \mathcal{X}$ represents the return of the position, i.e., $X(\mathbf{p})$ is the return of the position if scenario $\mathbf{p} \in \Omega$ happens. A measure of risk over $\mathcal{X}$ is a function that assigns a real-valued (or infinity) risk to each return function $X \in \mathcal{X}$.

DEFINITION 1. (Föllmera and Knispel 2013) A risk measure is a function $\rho : \mathcal{X} \to \mathbb{R} \cup \{+\infty\}$ with the following properties.

- Monotonicity: If $X \leq Y$ then $\rho(X) \geq \rho(Y)$.
- Translation Invariance: For any $m \in \mathbb{R}$, $\rho(X + m) = \rho(X) - m$.

These are natural properties that one might expect from a risk measure. If under any scenario, position $Y$ has a higher return than position $X$, i.e., $\forall \mathbf{p} \in \Omega, X(\mathbf{p}) \leq Y(\mathbf{p})$, then the risk associated to $X$ should be smaller than the risk of $Y$ (monotonicity property). Note that risk measure $\rho(X)$ is the monetary risk for holding the position $X$. Thus, if a real value $m$ is added to $X$, i.e., return of the position is increased by a constant value $m$ for all the scenarios, then the exposed risk should be reduced by $m$ (translation invariance).

Consider the state where we have served $h - 1$ impressions and the $h$-th impression arrives. This position can be described by $S = (w, s_1, s_2, \ldots, s_N)$; where $w$ is the revenue that has already been collected (sum of $r_{ji}$ for the allocated impressions), and $s_i$ is the remaining budget of the $i$-th advertiser. At each stage, we are uncertain about the future demands for available budgets $s_i$, $i = 1, 2, \ldots, N$. Or, equivalently, we are uncertain about the future price of each resource $p_i$. In the rest of the paper, we assume a measurable space $(\Omega, \mathcal{F})$ which is equipped with the Lebesgue measure $(P)$; $\Omega$ denotes the set of all possible prices $\Omega = \{\mathbf{p} | \mathbf{p} = (p_1, p_2, \ldots, p_N), 0 \leq p_i \leq p_{\max}\}$, and $\mathcal{Q}$ denotes the set of all probability measures on $(\Omega, \mathcal{F})$ which are absolutely continuous with respect to $P$.

Let $X_S : \Omega \to \mathbb{R}$ be the random variable that represents the value of the position $S$, i.e., for a given position $S = (w, s_1, s_2, \ldots, s_N)$, $X_S$ is a real-valued function on possible scenarios (prices) $\mathbf{p} \in \Omega$ that states the return of the position. If price $\mathbf{p}$ is realized for the resources, then value of the $i$-th resource would be $p_i s_i$. Thus, the total return would be

$$X_S(\mathbf{p}) = w + \sum_{i=1}^{N} p_i s_i.$$

We denote the set of all such functions by $\hat{\mathcal{X}} (\subset \mathcal{X})$.

$$\hat{\mathcal{X}} = \{X_S : \Omega \to \mathbb{R} \mid S = (w, s_1, s_2, \ldots, s_N), 0 \leq s_i \leq b_i, w \geq 0\},$$

As we alluded to before, the risk measure assigns a real number (or infinity) to each position $S$. From the definition of the risk measure, it should decrease linearly with the available cash (collected revenue $w$) and should also be decreasing with respect to the available budgets $s_i$. In this paper, for $X \in \mathcal{X}$, we study the class of the risk measures in the following form.

$$\rho(X) = \inf_{X \geq X_S, X_S \in \hat{\mathcal{X}}} \theta(X_S), \tag{5}$$

in which

$$\theta(X_S) = -w - \sum_{i=1}^{N} v_i(s_i). \tag{6}$$

In the above definition, $v_i(s_i) \in C^1$ is an increasing function with continuous derivative. For the consistency of the definition, we assume $\frac{\partial v_i}{\partial s}(s_i) \leq p_{\max}$; for more details see Lemma 6 in Appendix

A. We will show in Lemma 2, for $X_S \in \hat{\mathcal{X}}$, the risk measure in (5) coincides with the generating function $\theta(.)$. In other words, the risk measure $\rho(.)$ is the natural extension of $\theta(.)$ from $\hat{\mathcal{X}}$ to $\mathcal{X}$.

We start by showing that function $\rho(.)$, as it is defined in (5), satisfies conditions of a risk measure in Definition 1.

LEMMA 1. *Function $\rho : \Omega \to \mathbb{R} \cup \{+\infty\}$ defined in equation (5) is a risk measure.*

*Proof* We should establish the monotonicity and the transition invariant properties for the real valued function $\rho(.)$. The monotonicity property is clear from the definition (5). The transition invariance property can be shown as follows.

$$\rho(X+m) = \inf_{X \geq X_S - m, X_S \in \hat{\mathcal{X}}} \theta(X_S) \tag{7}$$

$$= \inf_{X \geq X_S, X_S \in \hat{\mathcal{X}}} \theta(X_S + m) \tag{8}$$

$$= \inf_{X \geq X_S, X_S \in \hat{\mathcal{X}}} \theta(X_S) - m. \tag{9}$$

Equation (9) is based on the observation that if $S = (w, s_1, s_2, \ldots, s_N)$ and $S' = (w + m, s_1, s_2, \ldots, s_N)$, then $X_S + m = X_{S'}$, and consequently $\theta(X_S + m) = \theta(X_S) - m$. Q.E.D.

While a risk measure should be defined over all possible return functions $X \in \mathcal{X}$, in the course of the allocation algorithm, we only encounter functionals $X_S \in \hat{\mathcal{X}}$. In the next lemma, we characterize risk measure $\rho(.)$ for return functions $X_S \in \hat{\mathcal{X}}$ which corresponds to a position $S$.

LEMMA 2. *In the impression allocation problem, for a given position $S = (w, s_1, s_2, \ldots, s_N)$, the risk measure (5) has a closed form representation as follows.*

$$\rho(X_S) = -w - \sum_{i=1}^{N} v_i(s_i).$$

*Proof* For $X_S \in \hat{\mathcal{X}}$,

$$\rho(X_S) = \inf_{X_S \geq X_{\bar{S}}, X_{\bar{S}} \in \hat{\mathcal{X}}} \theta(X_{\bar{S}}) \leq \theta(X_S) = -w - \sum_{i=1}^{N} v_i(s_i). \tag{10}$$

On the other hand, assume $X_{\bar{S}}$ is a position such that $X_S \geq X_{\bar{S}}$. By Lemma 6 in Appendix A

$$\theta(X_{\bar{S}}) \geq \theta(X_S). \tag{11}$$

Because equation (11) is valid for any $X_{\bar{S}}(\leq X_S)$,

$$\inf_{X_S \geq X_{\bar{S}}, X_{\bar{S}} \in \hat{\mathcal{X}}} \theta(X_{\bar{S}}) \geq \theta(X_S) = -w - \sum_{i=1}^{N} v_i(s_i). \tag{12}$$

Finally, the lemma follows from equations (10) and (12).     Q.E.D.

Having precisely defined the notion of risk in the display advertising problem, upon arrival of an impression, we can assign the impression to an advertiser such that the future risk is minimized. In the next section, we introduce an online resource allocation algorithm based on the risk minimization.

### 3.2. Risk Minimization Formulation

When an impression arrives, it can be assigned to any advertiser whose budget has not been exhausted yet. Each possible assignment would result in a different state, and subsequently a different risk. For an appropriately defined risk measure, we can allocate the impression such that the future risk is minimized.

Consider the position $S^{h-1} = (w^{h-1}, s_1^{h-1}, s_2^{h-1}, \ldots, s_N^{h-1})$ just before the $h$-th impression arrives. As detailed in the previous section, the risk associated to this position can be calculated by $\rho(X_{S^{h-1}})$. Now, if we assign the $h$-th impression to the $\hat{i}$-th advertiser, the new position would be

$$S^h = (w^{h-1} + r_{h\hat{i}}, s_1^{h-1}, s_2^{h-1}, \ldots, s_{\hat{i}-1}^{h-1}, s_{\hat{i}}^{h-1} - a_{h\hat{i}}, s_{\hat{i}+1}^{h-1}, \ldots, s_N^{h-1}),$$

and consequently the risk would be $\rho(X_{S^h})$.

To define $\mathcal{X}_{S^h}$ more precisely, let $b_i^{h-1}$ be the exhausted budget of the $i$-th advertiser right before the $h$-th impression arrives, i.e.,

$$b_i^{h-1} = \sum_{j=1}^{h-1} a_{ji} \hat{x}_{ji},$$

in which $\hat{x}_{ji}$ is one if impression $j$ was allocated to the $i$-th advertiser and zero otherwise. $\mathcal{X}_{S^h}$ is the set of random variables $X_{S^h}$ such that the corresponding position $S = (w, s_1, s_2, \ldots, s_n)$ satisfies the following feasibility conditions.

$$w = w^{h-1} + \sum_i r_{hi} x_{hi}$$

$$b_i^{h-1} + a_{hi}x_{hi} + s_i = b_i$$

$$\sum_i x_{hi} \leq 1$$

$$x_{hi} \geq 0, s_i \geq 0.$$

Therefore, in order to minimize the risk of the resulting position, the impression should be allocated to an advertiser by solving the following optimization:

$$\max \ \sum_i r_{hi}x_{hi} + \sum_i v_i(s_i) \tag{13}$$
$$\text{s.t.} \ \ b_i^{h-1} + a_{hi}x_{hi} + s_i = b_i$$
$$\sum_i x_{hi} \leq 1$$
$$x_{hi} \geq 0, s_i \geq 0;$$

where value function $v_i(.)$ of the remaining budget is described in (6) of the previous section.

Whenever an impression arrives, we solve minimization (13) to determine to which advertiser the impression should be allocated. As we discussed in Section 1, upon arrival of the impression, it should be assigned to an advertiser almost immediately. Later in Section 5, we show the optimal solution in (13) has a closed form solution and it can be calculated efficiently. Furthermore, note that each value function $v_i(.)$ induces an allocation algorithm. In the next section, we use the robust representation of the risk measure to relate the auctioneer's belief of the future demand to the choice of the budget value function $v_i(.)$. In other words, the auctioneer (ad network) would choose a budget value function based on his or her belief regarding the uncertainty in the future.

## 4. Convex Risk Measure and Robust Representation

Introduction of a convex risk measure in the context of ad allocation is the key step for building a link between the auctioneer's belief (of the future demands) and the budget value function $v_i(.)$ in equation (5). Later in Section 5, we show a belief in the future distribution of impressions is translated to a price updating rule in the allocation algorithm.

A convex risk measure, as defined below, has a robust (dual) representation that provides valuable insights to the corresponding risk measure. In this section, we first define the notion of convexity for a risk measure, and then the robust representation of a risk measure is introduced. Next, we find this representation for the class of risk measures that are introduced in Section 3.

DEFINITION 2. (Föllmera and Knispel 2013) A risk measure $\rho(.)$ is called convex if it satisfies the following convexity condition.

$$\forall \lambda \in [0,1], \ \rho(\lambda X + (1-\lambda)Y) \leq \lambda \rho(X) + (1-\lambda)\rho(Y).$$

The convexity condition intuitively means that the risk of two positions combined should be less than combination of individual risks, i.e., diversification decreases the risk.

LEMMA 3. *If the budget value functions $v_i(.)$ are concave and increasing then the class of risk measures defined in (5) are convex.*

*Proof* First, note that since $v_i(.)$ are concave and increasing then, by definition, $\theta(.)$ is a convex and decreasing function.

Let $X, Y \in \mathcal{X}$ be two return functions. By definition of the risk measure $\rho(.)$, for each $\epsilon_1, \epsilon_2 > 0$, there exist $X_{S_1}, X_{S_2} \in \hat{\mathcal{X}}$ such that

$$X \geq X_{S_1}, Y \geq X_{S_2}$$

$$\rho(X) = \inf_{X \geq X_S, X_S \in \hat{\mathcal{X}}} \theta(X_S) \geq \theta(X_{S_1}) - \epsilon_1$$

$$\rho(Y) = \inf_{Y \geq X_S, X_S \in \hat{\mathcal{X}}} \theta(X_S) \geq \theta(X_{S_2}) - \epsilon_2.$$

Therefore,

$$\begin{aligned}
&\lambda \rho(X) + (1-\lambda)\rho(Y) \\
&= \lambda \inf_{X \geq X_S, X_S \in \hat{\mathcal{X}}} \theta(X_S) + (1-\lambda) \inf_{Y \geq X_S, X_S \in \hat{\mathcal{X}}} \theta(X_S) \\
&\geq \lambda \theta(X_{S_1}) - \epsilon_1 + (1-\lambda)\theta(X_{S_2}) - \epsilon_2 \\
&\geq \theta(\lambda X_{S_1} + (1-\lambda) X_{S_2}) - \epsilon_1 - \epsilon_2 \\
&\geq \inf_{\lambda X + (1-\lambda) Y \geq X_S, X_S \in \hat{\mathcal{X}}} \theta(X_S) - \epsilon_1 - \epsilon_2 \\
&= \rho(\lambda X + (1-\lambda)Y) - \epsilon_1 - \epsilon_2.
\end{aligned} \quad (14)$$

In the above equations, line 3 to 4 is due to the convexity of $\theta(.)$, and line 4 to 5 follows from the fact that $\lambda X_{S_1} + (1-\lambda) X_{S_2} \in \hat{\mathcal{X}}$. Since equation (14) is valid for every $\epsilon_1, \epsilon_2 > 0$, the convexity property follows.   Q.E.D.

In Lemma 3, we showed that the risk measure $\rho(.)$ is convex. Next, for a given family of probability measures $\mathcal{Q}$ on possible scenarios $\Omega$, we prove that the risk measure $\rho(.)$ has a robust (dual) representation (Föllmera and Schied 2002, Föllmera and Knispel 2013). Intuitively, the risk measure of a return function $X$ would be interpreted as the maximum loss under any probability measure on $(\Omega, \mathcal{F})$, however, there would be a penalty $\alpha(Q)$ associated to each probability measures $Q \in \mathcal{Q}$. Penalty function $\alpha(Q)$ represents our belief toward the probability measure. Equivalently, the risk measure can be defined based on the penalty function on probability measures. Robust representation is one of the most powerful properties of a convex risk measure and it provides deep insights to the underlying risk measure.

LEMMA 4. *Let $\mathcal{Q}$ be a set of probability measures on $(\Omega, \mathcal{F})$. If value functions $v_i(.)$ are continuous, then risk measure $\rho(.)$ in (5) can be represented by a probability measure penalty function (in short, probability penalty) $\alpha(.)$.*

$$\rho(X) = \sup_{Q \in \mathcal{Q}} \Big( E_Q[-X] - \alpha(Q) \Big), \tag{15}$$

*where*

$$\alpha(Q) = \sup_{X \in \mathcal{X}} \Big( E_Q[-X] - \rho(X) \Big). \tag{16}$$

*Proof* A convex risk measure has a robust representation (15) if it is almost surely continuous from above (see Theorem 4.30 in Föllmera and Schied (2004)). That is,

$$X_n \searrow X \implies \rho(X_n) \nearrow \rho(X).$$

Details of the proof can be found in Appendix B.   Q.E.D.

In the robust representation, risk of a random variable $X$ is equal to the maximum expected loss under the probability measure $Q$, $E_Q[-X]$, with penalty function $\alpha(Q)$ for choosing $Q$.

In the next lemma, we characterize the penalty function $\alpha(Q)$ for the class of risk measures $\rho(.)$ defined in (5). Note that, in equation (16), the supremum is taken over all $X \in \mathcal{X}$ which includes all real-valued bounded functions on $\Omega$. While for a generic risk measure and function space $\mathcal{X}$, we cannot further simplify the representation of $\alpha(.)$, in the impression allocation problem, the penalty function $\alpha(.)$ has a much simpler form.

LEMMA 5. *Let $Q \in \mathcal{Q}$ be a probability measure on $(\Omega, \mathcal{F})$. Then, probability penalty function $\alpha(Q)$ can be restated by taking suppremum over $\hat{\mathcal{X}}$ as follows.*

$$\alpha(Q) = \sup_{X_S \in \hat{\mathcal{X}}} \sum_{i=1}^{N} \left( v_i(s_i) - p_i^* s_i \right) \tag{17}$$

*where*

$$p_i^* = \sum_{\mathbf{p} \in \mathcal{P}} Q(\mathbf{p}) p_i.$$

*Proof* Appendix C.

In the resource allocation problem, $Q$ is a probability measure on the possible future resource prices (or equivalently, the demand for the remaining budgets $s_i$). Thus, $p_i^*$ is the expected dual price for the $i$-th advertiser assuming probability measure $Q$ for future demand.

Note that each probability penalty function $\alpha(Q)$ is associated to a budget value function $v_i(.)$. In what follows, we consider four budget value functions and find their representations in term of the probability penalty function.

**Zero Value Function:** If we choose the budget value function to be zero, $v_i(s_i) = 0$ for all $i$, that is, future opportunities for spending the budgets are ignored, and only the immediate revenue, $w$, is valued.

$$\rho(X_S) = -w. \tag{18}$$

While budget constraints are not violated, each impression is allocated to the advertiser with the highest value (bid). Choosing zero value function reduces the allocation algorithm to the greedy allocation and $\alpha(Q) = 0$. When $\alpha(Q)$ is zero, we are not penalizing any probability measure, and we consider the worse case for evaluating the future opportunities. The performance of greedy allocation is widely studied in the literature (Mahdian 2011, Karande et al. 2011, Mehta et al. 2005).

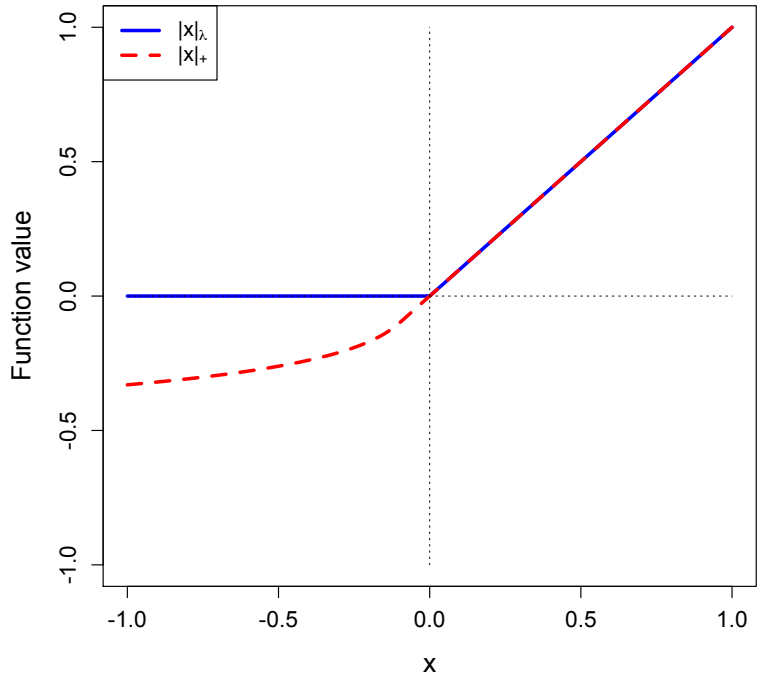

**Figure 2** Comparing two functions : $|x|_+$ (used in linear budget value function) and $|x|_\lambda$ (used in logarithmic budget value function).

As mentioned earlier, advertisers often require their budget to be spent smoothly through the flight of the campaign. The greedy algorithm not only doesn't compensate for the future demands, it would exhaust budgets of advertisers non-uniformly. Typically, at the beginning, advertisers with high valuations (bids) would get most of the impressions until they are out of budget. Then, impressions would be allocated to advertisers with lower bids.

**Linear Value Function:** Let $p_i^\epsilon$ be an estimate of resource prices in the future. For example, $p_i^\epsilon$ can be the dual price of the $i$-advertiser's budget from the sample linear program based on previously observed impressions. The budget value function can be defined as $v_i(s_i) = p_i^\epsilon s_i$ for all $i$, and, therefore,

$$\rho(X_S) = -w - \sum_{i=1}^{N} p_i^\epsilon s_i. \tag{19}$$

By Lemma 5, $\alpha(Q)$ is

$$\begin{aligned}
\alpha(Q) &= \sup_{X \in \mathcal{X}} \sum_{i=1}^{N} \left( v_i(s_i) - p_i^* s_i \right) \\
&= \sup_{X \in \mathcal{X}} \sum_{i=1}^{N} \left( p_i^\epsilon s_i - p_i^* s_i \right) \\
&= \|\mathbf{p}^\epsilon.\mathbf{b} - \mathbf{p}^*.\mathbf{b}\|_+.
\end{aligned} \tag{20}$$

Note that dot operation (.) in both $\mathbf{p}^\epsilon.\mathbf{b}$ and $\mathbf{p}^*.\mathbf{b}$ is component-wise vector multiplication, and $\|.\|_+$ is the norm-one of the positive components.

$$\|\mathbf{x}\| = \sum |x_i|_+,$$

where

$$|x|_+ = \begin{cases} x \text{ if } x \geq 0; \\ 0 \text{ if } x < 0. \end{cases}$$

Intuitively, if the estimated price $p_i^\epsilon$ is less than the expected $i$-th price from the probability measure, $p_i^*$, we have underestimated the value of the resource $s_i$, i.e., if probability measure $Q$ explains the future demand, the final revenue would be higher than $p_i^\epsilon s_i$. Hence, expected revenue $E_Q[-X]$ is not penalized for this resource. On the other hand, if the estimated price $p_i^\epsilon$ is greater than the expected $i$-th price, $p_i^*$, then our estimation for the collected revenue from the $i$-th resource is larger than what we would collect under $Q$. Thus, $E_Q[-X]$ is penalized linearly.

The above two budget value functions are static in the sense that, while budget constraints are satisfied, the allocation algorithm doesn't depend on the remaining budgets of advertisers. In other words, a small error in the initial estimation of the dual prices can propagate and lead to unbalanced budget consumption. In the next two algorithms, we choose budget value functions $v_i(s)$ such that the allocation algorithm is dynamically adjusted, i.e., if the budget of a particular advertiser is being exhausted faster than it should be spent, the algorithm reacts to the situation by reducing the number of impressions allocated to the advertiser.

**Log Value Function:** Again, assume $p_i^\epsilon$ is an estimate of resource prices. In the case of the linear value function, the budget of the $i$-th advertiser would be spent much faster when $p_i^\epsilon$ is smaller than the realized price of the $i$-th resource. Our goal is to dynamically adjust the allocation algorithm in case of overspending, and this can be achieved by appending the budget value function with a logarithmic term.

$$v_i(s_i) = p_i^\epsilon s_i + \gamma \log\left(\frac{s_i}{b_i}\right), \ i=1,\ldots,N. \tag{21}$$

In this value function, $\gamma$ is a regularization factor that can be used to tune the algorithm. In practice, we can add a small constant value to the term inside $\log(.)$ in order to bound the derivative of $v_i(.)$ and enforce $\frac{\partial v_i}{\partial s}(s) \leq p_{\max}$. For the log value function above, the probability penalty function $\alpha(Q)$ is (see Appendix D)

$$\alpha(Q) = \|\mathbf{p}^\epsilon.\mathbf{b} - \mathbf{p}^*.\mathbf{b}\|_\gamma,$$

where

$$\|\mathbf{x}\|_\gamma = \sum |x_i|_\gamma,$$

and

$$|x|_\gamma = \begin{cases} x & \text{if } x \geq -\gamma; \\ -\gamma + \gamma \log\left(\frac{\gamma}{-x}\right) & \text{if } x < -\gamma. \end{cases} \tag{22}$$

The log value function is similar to the linear value function when $p_i^\epsilon - p_i^* \geq 0$; the probability measure is linearly penalized by $p_i^\epsilon - p_i^*$. However, we negatively penalize (reward) probability measures in which $p_i^*$ is greater than the estimated price $p_i^\epsilon$. In other words, we are more confident about the probability measures implying that our initial estimates for the prices $p_i^\epsilon$ are smaller than what would be realized. Figure 2 compares two functions $\|\mathbf{x}\|_+$ (used in linear value function) and $\|\mathbf{x}\|_\lambda$ for $\lambda = 0.1$. While the functions coincide for $x \geq 0$, they behave differently in the negative regime.

This resulting algorithm heavily penalizes advertisers with small remaining budgets $s_i$. As the available budget $s_i$ approaches zero, the value function $v_i(s_i)$ approaches negative infinity. Thus, we tend to avoid allocating impressions to advertisers with small remaining budgets. In other words, the algorithm focuses on conserving budgets and allocating impressions to advertisers with higher available resources.

**Exponential Value Function:** Another approach is to use an exponential function to dynamically update the estimated dual prices. When the algorithm is overspending the budget of a particular advertiser, the initial estimate of the dual price is increased and, conversely, the initial estimate of the dual price is decreased in the case of underspending. The exponential budget value function is defined as follows.

$$v_i(s_i) = p_i^\epsilon \left( -\frac{b_i}{\kappa} \right) e^{k\left( \frac{b_i - s_i}{b_i} - \frac{h}{N} \right)}, \ i = 1, \ldots, N, \tag{23}$$

where $\kappa > 0$ is a regularizing parameter. In the exponential value function, $(b_i - s_i)/b_i$ is the fraction of the exhausted budget after the allocation and $h/N$ is the fraction of allocated impressions. For $\kappa > 0$, if the fraction of used budget, $(b_i - s_i)/b_i$, is more than the fraction of allocated impressions, $h/N$, the value function $v_i(.)$ would be small. Consequently, less impressions would be allocated to these advertisers.

For the exponential budget value function, the penalty for the probability measure $Q$ is a function of two quantities: (a) KL-divergence between $\mathbf{p}^\epsilon.\mathbf{b}$ and $\mathbf{p}^*.\mathbf{b}$, and (b) fraction of total unserved impressions $1 - \frac{h}{N}$ (see Appendix E).

$$\alpha(Q) = \left( \sum_{i=1}^{N} p_i^* b_i \right) \left( -\left(1 - \frac{h}{N}\right) + \frac{1}{\kappa} \left( \mathcal{D}_{KL}\left(\mathbf{p}^\epsilon.\mathbf{b} \| \mathbf{p}^*.\mathbf{b}\right) - 1 \right) \right).$$

Note that, in computing KL-divergence, we require $\mathbf{p}^\epsilon.\mathbf{b}$ and $\mathbf{p}^*.\mathbf{b}$ to be normalized. As we will show in Section 6, among all the introduced allocation algorithms, this algorithm has the best balance between maximizing the immediate revenue and conserving budgets for future impressions.

# 5. Updating Dual Prices with Risk Minimization

In Section 3, we described the risk minimization framework for the impression allocation problem. In practice, an adserver receives tens of thousands of impression requests per second, and the assignment should happen almost immediately, typically in less than 5 milliseconds. Thus, the allocation algorithm should be computationally inexpensive. In this section, we show that the risk minimization in (13) has a closed form solution, and each risk measure is translated to a simple dual price updating rule. As the result, impressions can be allocated in realtime without heavy computations.

Recall that when the $h$-th impression is received, we assign it to an advertiser to minimize the risk measure in (6), or equivalently, maximize instantaneous revenue, $\sum_i r_{ji} x_{ji}$, plus a value function, $\sum_i v_i(s_i)$, that controls overspending of advertisers. Here, variable $s_i$ represents the remaining budget of the $i$-th advertiser after assigning the $h$-th impression. Thus, we precisely solve optimization problem (13) upon arrival of the $h$-th impression.

The solution to optimization problem (13) can be characterized by its Lagrange function.

$$L(x_{hi}, s_i) = \sum_i r_{hi} x_{hi} + \sum_i \nu_i(b_i - b_i^{h-1} - a_{hi} x_{hi}) + \lambda(\sum_i x_{hi} - 1).$$

In any optimal solution, we should have $\frac{\partial L}{\partial x_{hi}} \geq 0$.

$$\frac{\partial L}{\partial x_{hi}} = r_{hi} + \nu_i'(b_i - b_i^{h-1} - a_{hi} x_{hi})) \times (-a_{hi}) + \lambda.$$

Assume $\mathbf{x}_h^*$ is an optimal solution. If $x_{hi}^* > 0$ then $\frac{\partial L}{\partial x_{hi}} = 0$, and therefore,

$$i \in \mathrm{argmax}_k \ r_{hk} + \nu_i'(b_k - b_k^{h-1} - a_{hk} x_{hk}) \times (-a_{hk}).$$

If $\max_k r_{hk} + \nu_i'(b_k - b_k^{h-1} - a_{hk} x_{hk}) \times (-a_{hk}) \leq 0$ then we will not allocate the impression to any advertiser and $\hat{x}_{hi} = 0$ for all $i$. However, in practice, ad networks often have a *house* campaign with an infinite budget and very small revenues $(r_{ji})$. Thus, $\max_k r_{hk} + \nu_i'(b_k - b_k^{h-1} - a_{hk} x_{hk}) \times (-a_{hk})$ is always positive in practice.

Next, we find closed form solutions for the four budget value functions that are discussed in Section 4. In all cases, the allocation algorithm is equivalent to a simple price updating scheme.

**Zero Value Function (greedy):** In this case we choose $v_i(s_i) = 0$ for all $i$, so that the allocation algorithm would reduce to the greedy algorithm. Upon arrival of the $h$-th impression request, it would be allocated to the $i^*$-th advertiser according to (break ties arbitrarily)

$$i^* \in \mathrm{argmax}_i \quad r_{hi} \tag{24}$$

$$\text{s.t.} \quad b_i^{h-1} + a_{hi} \le b_i, \ i = 1, 2, \ldots, M. \tag{25}$$

Set $\hat{x}_{hi^*} = 1$ and $\hat{x}_{hi} = 0$ for $i \ne i^*$.

**Linear Value Function (fixed dual):** This is the case we choose $v_i(s_i) = p_i^\epsilon s_i$ for all $i$, so that the allocation algorithm would be equivalent to the dual-based algorithm proposed in Agrawal et al. (2009).

$$i^* \in \mathrm{argmax}_i \quad r_{hi} - p_i^\epsilon a_{hi} \tag{26}$$

$$\text{s.t.} \quad b_i^{h-1} + a_{hi} \le b_i, \ i = 1, 2, \ldots, N. \tag{27}$$

Effectively, the dual price of the $i$th advertiser is fixed.

**Log Value Function (log):** This is the case we choose $v_i(s_i) = p_i^\epsilon s_i + \gamma \log\left(\frac{s_i}{b_i}\right)$. Since $v_i'(s_i) = p_i^\epsilon + \frac{\gamma}{s_i}$, the allocation algorithm would be

$$i^* \in \mathrm{argmax}_i \quad r_{hi} - p_i^\epsilon a_{hi} - \gamma \frac{a_{hi}}{s_i} \tag{28}$$

$$\text{s.t.} \quad b_i^{h-1} + a_{hi} + s_i = b_i, \ i = 1, 2, \ldots, M. \tag{29}$$

$$s_i \ge 0. \tag{30}$$

Effectively, the dual price of the $i$th advertiser is dynamically adjusted to $p_i^\epsilon + \frac{\gamma}{s_i}$.

**Exponential Value Function (exponential):** This is the case we choose

$$v_i(s_i) = p_i^\epsilon \left(-\frac{b_i}{\kappa}\right) e^{k\left(\frac{b_i - s_i}{b_i} - \frac{h}{N}\right)},$$

so that the allocation algorithm would be

$$i^* \in \text{argmax}_i \quad r_{hi} - p_i^\epsilon a_{hi} e^{\kappa(\frac{b_i - s_i}{b_i} - \frac{h}{N})} \tag{31}$$

$$\text{s.t.} \quad b_i^{h-1} + a_{hi} + s_i = b_i, \ i = 1, 2, \ldots, M.$$

$$s_i \geq 0.$$

Effectively, the dual price of the $i$th advertiser is dynamically adjusted to $p_i^\epsilon e^{\kappa(\frac{b_i - s_i}{b_i} - \frac{h}{N})}$.

## 6. Experimental Results

While online resource allocation is widely studied under various theoretical assumptions, the performance of these proposed methods is not comprehensively analyzed with real data. We used allocation algorithms discussed in Section 5 for the revenue maximization in the display advertising where billions of impressions are served every day. First, in Section 6.1 we describe the data structure and characteristics. In Section 6.2 we discuss dual prices estimation using the previously observed data. Next, the performance of various allocation algorithms is presented in Section 6.3. As mentioned before, usually advertisers enforce some restrictions on the number of served ads to a particular user (frequency capping). In Section 6.4, we extend the framework to accommodate ad serving under frequency capping constraints.

### 6.1. Data and Experiment Setting

We used AOL/Advertising.com impression logs to evaluate the performance of the proposed price updating algorithms. The data set includes about 200 million impressions that were served on October 11, 2013 between 10 am and 12 am. We restrict this study to 700 advertisers (campaigns) with limited budgets. Impressions were divided into two subsets:

$T_1$ : Impressions which were served between 10 am- 11 am.

$T_2$ : Impressions which were served between 11 am- 12 am.

Impressions in $T_1$ were used for learning initial dual prices in equation (2), and impressions in time period $T_2$ were used for evaluation.

Note that there are around 100 million impressions in $T_1$, and that it is not possible to incorporate all impressions in the LP solver (equation (2)). If we have $M$ impressions, then the number of constraints in LP would be $M + N$ (N=700) and the number of variables would be $MN$. With millions of impressions, the LP size is too large for the solver to handle. Thus, we sample impressions in $T1$ with the rate of $\delta$ which generates a smaller LP.

**Non-zero Bids of Each Advertiser:** When an impression request arrives, many advertisers submit zero value for the impression. This can be due to a number of factors such as geographic targeting, or audience targeting. For example, a campaign might be promoting a particular product to female users in the US with an income of more than $100k. Thus, the advertiser assigns zero value for most of the impressions as he/she only values the narrow subset of impressions coming from the target group.

Figure 3 is the distribution of advertisers based on the submitted non-zero bids for about 910k impressions. As one can see, many advertisers (more than 170 out of 700) submit less than 50k bids. These advertisers are targeting a very specific group of Internet users.

**Non-zero Bids of Each Impression:** Narrow targeting can also be observed by focusing on the number of non-zero bids for each impression. Figure 4 shows frequency of impressions based on submitted bids (non-zero). As one can see in the figure, there is no impression that is valued by more than 450 advertisers. Indeed, almost all of the impressions receive less than 400 (out of 700) non-zero bids, and almost half of the impressions receive less than 200 bids.

Figures 3 and 4 show that the initial LP problem is sparse. Some advertisers are only interested in a very small fraction of the impressions. Generally, it is much harder to allocate the budget of these narrowly targeted advertisers efficiently.

**LP solver:** We used MOSEK as the LP solver in phase 1 for estimating dual prices from impressions in $T_1$. Constraining the running time to less than 1.5 hours, and using a machine with an Intel Core i7 processor and 8G of RAM, we could solve instances of the problem with up to 52k impressions.

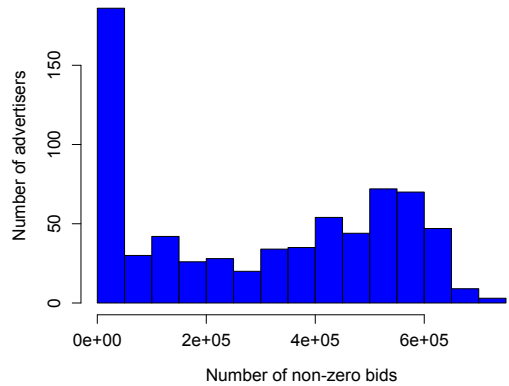

**Figure 3** Frequency of advertisers based on non-zero bids that they submit for about 910k impression requests.

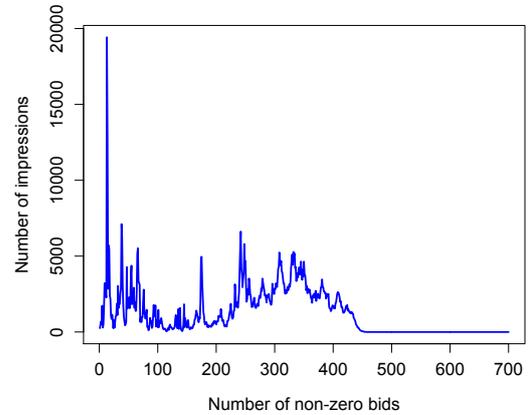

**Figure 4** Frequency of impressions based on the number of advertisers who submit non-zero bids for them.

### 6.2. Estimating Dual Prices in the Initial Phase

Since initial dual prices $p_i$, $i = 1, 2, \ldots, N$, are estimated from a subset of impressions in $T_1$, estimated dual prices depend on a randomly selected subset of impressions. As you may expect, the stability of estimated dual prices is a key factor in the performance of the allocation algorithms. Here, we show how dual prices are impacted as we add more impressions to the LP in equation (2).

In order to study the volatility of dual prices, the set of impressions in $T_1$ is permuted randomly, and then LP is solved using the first $M_l$ impressions in the permuted set. Figure 5 compares dual prices that are computed using a different number of impressions. In the top left figure, we compare dual prices computed using $M_l = 100$ impressions versus dual prices computed using $M_l = 200$ impressions. Note that the set of $M_l = 200$ impressions is obtained by appending 100 randomly-selected impressions to the initial 100 impressions. Thus, the two sets of impressions are not independent. The size of each point is proportional to $\log K_i$, where $K_i$ is the number of non-zero bids of the $i$-th advertiser. Smaller points correspond to advertisers with a small number of non-zero bids, and it is typically harder to find their dual prices.

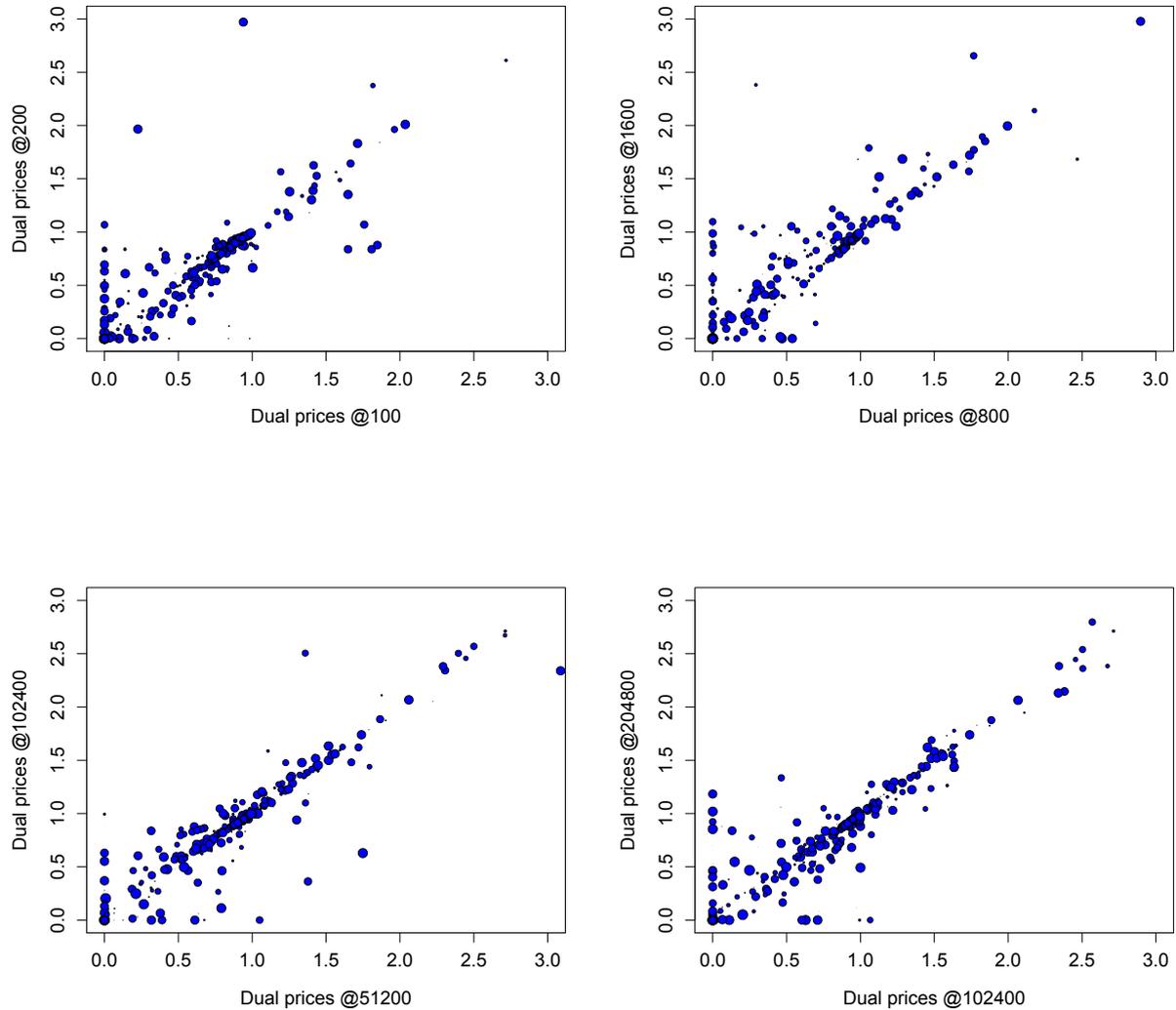

**Figure 5** Dual prices from different sets of impressions. Top left: dual prices after 100 impressions vs 200 impressions. Top right: dual prices after 800 impressions vs 1600 impressions. Bottom left: dual prices after 6400 impressions vs 12800 impressions. Bottom right: dual prices after 204800 impressions vs 409600 impressions. Note that the size of each point is proportional to log(number of non-zero bids of the advertiser).

The top right figure represents dual prices at $M_l = 800$ versus $M_l = 1600$, and down left figure represents dual prices at $M_l = 6,400$ versus $M_l = 12,800$. Finally, the down left plot in Figure 5 compares dual prices at $M_l = 204,800$ and $M_l = 409,600$. It is important to note that using MOSEK, we could only include at most 52k impressions in the LP (to be solved in less than 1.5 hours). In the case of $M_l > 52k$, we just used the last 52k impressions. As a result, in contrast to the first three figures, impression sets in the last figure are disjointed.

As the number of impressions are increased in the LP, most of the dual prices converge. That is, they have the same value for the two different impression sets. However, for some advertisers, dual prices can be significantly different in two samples, i.e., the representing points in Figure 5 are far from the identity line. Even in the last plot (bottom left), which includes the most number of impressions from $T_1$, there are many advertisers whose dual prices are zero in one sample and significant in the other sample. This observation indicates that using static dual prices for time $T_2$ can result in an unbalanced budget consumption for some advertisers. As one might expect, dynamic adjustment of dual prices in log value and exponential value can mitigate the effect of unstable dual prices. We will discuss this effect in more detail in the next section.

### 6.3. Performance of Proposed Allocation Algorithms

Ad networks usually have two criteria for the performance of an allocation algorithm that are closely connected. The first criterion is the total collected revenue as discussed earlier. Indeed, we only used revenue criterion to formulate and study the problem of display advertising allocation. In addition to the revenue performance, advertisers require their budget to be spent evenly. Typically, it is not desirable to exhaust the budget of an advertiser long before the end of the day. Using real data, we demonstrate how various value functions perform under these metrics.

We formulated the display advertising allocation as a linear programming problem with the total revenue as the objective function (equation (1)). Figure 6 shows the total collected revenue of four allocation algorithms as they progress. The greedy algorithm, as the name suggests, is greedy in collecting revenue at the beginning; however, soon fails to keep pace and total revenue starts to saturate. Using fixed duals for allocation, based on equation (19), would regularize bids and result in a higher revenue. Note that dual prices are estimated from the previous data points ($T_1$) and that they don't totally fit the new data set $T_2$. The next allocation method is based on log value function (21); the total revenue is less than fixed duals and more than greedy. As we will discuss later, log value induces an overly conservative allocation method, and does not spend the advertiser's total budget. Finally, we have exponential value function which has the highest revenue. This method uses estimated dual prices and keeps updating them according to the served impressions. Table 1 contains performance metrics of four algorithms. The total revenue of the exponential allocation is 11.6% better than the greedy algorithm.

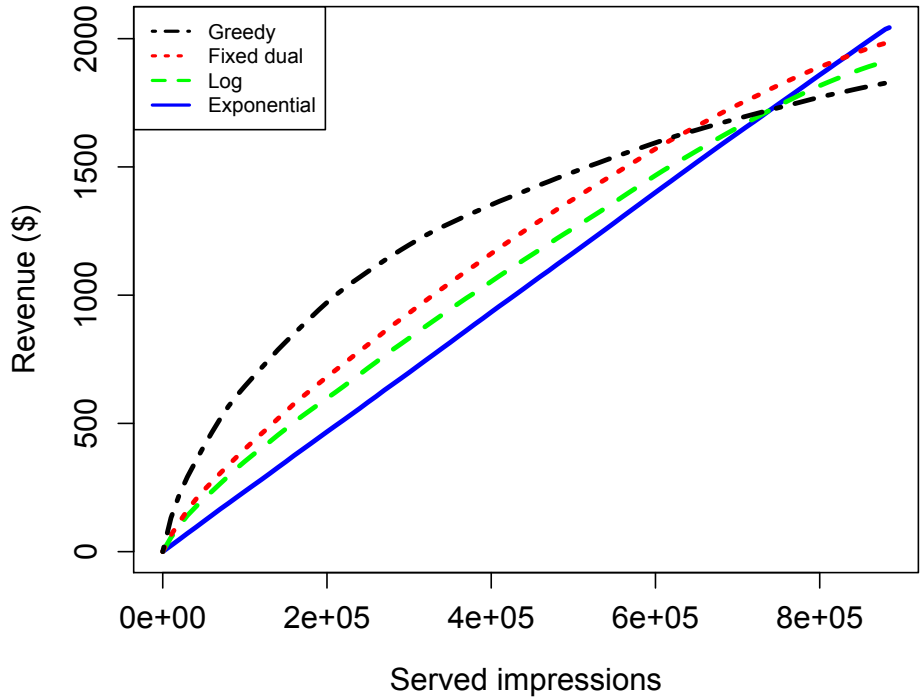

**Figure 6**     The total revenue of four allocation algorithms. Greedy algorithm generates high revenue at the beginning and saturates very fast. While performance of fixed dual prices is reasonable, log value is over-conservative in preserving resources. Allocation corresponding to the exponential value function achieves the highest revenue and smooth spending.

In addition to the revenue, balanced budget delivery is regarded critical for advertisers and ad networks. Advertisers prefer to reach their audiences regularly during the flight of a campaign rather than run out of budget early before the end. The number of campaigns which are out of budget (oob) for the four discussed allocation algorithms are shown in Figure 7. Greedy allocation assigns each impression to the advertiser with the highest immediate revenue. Therefore, as expected, the greedy algorithm aggressively exhausts budgets of campaigns with higher revenue. As you can see in Table 1, by mid flight, 366 campaigns would be out of budget (oob) using greedy allocation.

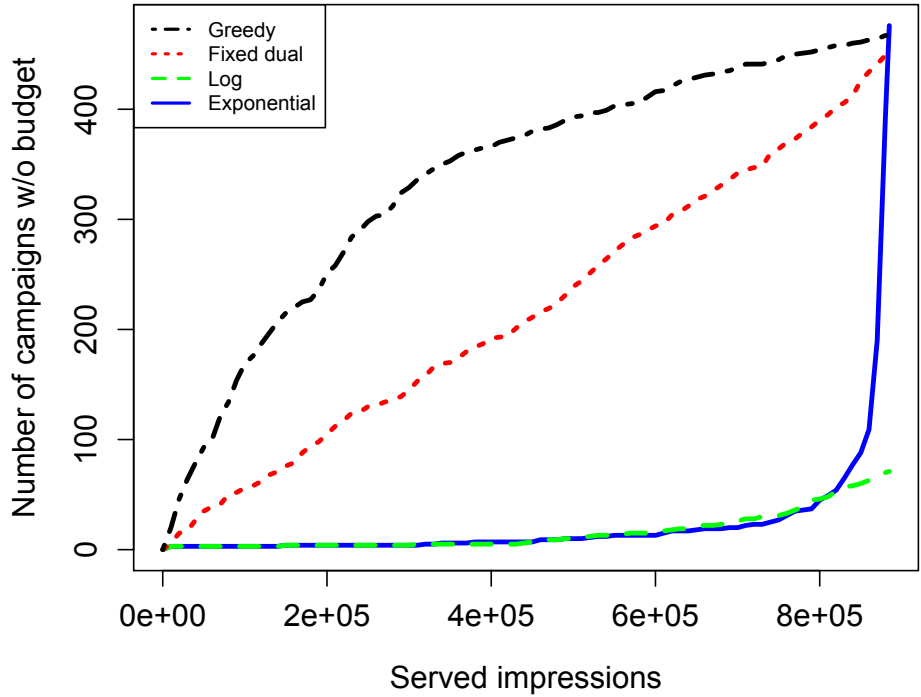

**Figure 7** The total number of campaigns with exhausted budget vs the number of served impressions. In both greedy and fixed dual, advertisers start to exhaust their budget from the beginning of the experiment. In contrast, log value and exponential value allocations keep advertisers in budget until the very end.

Discounting bids of advertisers with dual prices not only improves the total allocation revenue but it also results in an smoother delivery. If an advertiser with a limited budget has a high value for most impressions, then its associated dual price would be high and thus regularized bids, $r_{ji} - p_i a_{ji}$, would be much smaller. Therefore, the advertiser would get less impressions and its budget is conserved during the flight. Using dual prices, the number of campaigns with no budget (oob) is 192 in mid flight compared to 366 oob campaigns in greedy (Table 1).

**Table 1**    Performance of allocation algorithms.

| Allocation algorithm | Total Revenue | Improvement over greedy | Mid flight oob | Final oob |
|---|---|---|---|---|
| Greedy | $1829.94 | - | 366 | 467 |
| Fixed dual | $1986.67 | 8.5% | 192 | 452 |
| Log | $1915.72 | 4.6% | 5 | 71 |
| Exponential | $2043.21 | 11.6% | 7 | 476 |

The last two allocation algorithms, log value and exponential value, ensure a very smooth delivery. Log value function, in equation (21), highly penalizes bids of advertisers with small remaining budgets. Thus, as it is clear from Figure 7 and Table 1, most of the advertisers stay in budget during the flight. This results in an overly conservative allocation algorithm. Even at the end of the experiment, there are only 71 campaigns with no budget, which is not desirable in terms of the total revenue.

The exponential value function has the best balance between attaining the highest revenue and smooth delivery. As one can see in Table 1, while it has the highest total revenue (11.6% improvement over greedy), most of the advertisers are in budget during the flight. It is only at the end of the flight that exponential value starts to exhaust budgets of advertisers, which is a desirable behavior.

To test robustness of the algorithms further, we applied four allocation algorithms to the stream of the impressions in the reverse order as well. Impressions are fed into the algorithms from last to first, and the results are shown in Table 2. Performance metrics are very similar to what we observed in the natural order of impressions.

### 6.4. Enforcing Frequency Caps (Fcap)

Showing a particular ad to one user multiple times might be less effective than showing the same number of ads to many distinct users. Thus, advertisers typically enforce some frequency caps for showing a specific ad to one user. For example, they may require that each user is served a particular ad at most two times per 24 hours. In this section, we extend the proposed methods to be able to handle frequency capping.

**Table 2** Performance of allocation algorithms for the reverse impression stream.

| Allocation algorithm | Total Revenue | Improvement over greedy | Mid flight oob | Final oob |
|---|---|---|---|---|
| Greedy | $1828.34 | - | 380 | 467 |
| Fixed dual | $1988.09 | 8.7% | 214 | 455 |
| Log | $1916.11 | 4.8% | 8 | 68 |
| Exponential | $2044.02 | 11.7% | 8 | 473 |

Frequency capping can be enforced in various ways. In most cases, advertisers restrict the ad network to show at most $f$ ads in each $h$ hours. For simplicity, assume, for the given time period that we could serve at most $f_i$ ads from the $i$-th advertiser to each user. If $\mathcal{I}_u$ is the set of impressions from user $u$, then we can solve the following the linear optimization problem to estimate dual prices in phase 1.

$$\max \sum_{j=1}^{\epsilon M} \sum_{i=1}^{N} r_{ji} x_{ji} \quad (32)$$

$$\text{s.t.} \sum_{j} a_{ji} x_{ji} \leq b_i^\epsilon$$

$$\sum_{j \in \mathcal{I}_u} x_{ji} \leq f_i \ \forall u \in 1, 2 \ldots, U \ \forall i \in 1, 2, \ldots, N \quad (33)$$

$$\sum_{i} x_{ji} \leq 1, x_{ji} \geq 0$$

However, adding frequency capping constraints (33) to the LP formulation is not practical. If we have $U$ users, we should add one constraint per user per advertiser. Thus, we should add $UN$ new constraints to the original LP. The number of users, $U$, is usually very large and proportional to the number of impressions $M$. Current commercial LP solvers wouldn't be able to solve it in a reasonable time (e.g. less than 1.5 hours in our application) with tens of thousands of impressions.

In order to reduce the computation complexity of the LP in (32), we relax the LP (32) to reduce the number of constraints while frequency caps constraints are approximately enforced. For advertiser $i$, let $\mathcal{P}_i$ be a partition over impressions with non-zero bids. Instead of enforcing a user-based frequency cap, we establish a frequency limit on each partition. In other words, for each set $P \in \mathcal{P}_i$, we compute an upper bound, $f(i, P)$, for the number of impressions that can be served to the particular set $P$.

$$f(i,P) = \max \ |S| \qquad (34)$$
$$\text{s.t.} \ \ S \subset P,$$
$$\sum_{j \in \mathcal{I}_u \cap S} x_{ji} \leq f_i \ \ u = 1, 2, \ldots, U,$$

where $|.|$ is the size of set $S$. Note that in this optimization problem, we want to find the maximum number of impressions from set $P$ that can be assigned to an advertiser while frequency caps are satisfied. Since constraints are on individual users, the optimal value, $f(i, P)$, can be computed in $O(|P|)$.

Finally, we can replace the user based frequency capping constraints in problem (32) with partition based frequency caps.

$$\max \ \sum_{j=1}^{\epsilon M} \sum_i r_{ji} x_{ji} \qquad (35)$$
$$\text{s.t.} \ \sum_j a_{ji} x_{ji} \leq b_i^\epsilon,$$
$$\sum_{j \in P_i} x_{ji} \leq f(i, P_i) \ \forall P_i \in \mathcal{P}_i \ \forall i \in 1, 2, \ldots, N,$$
$$\sum_i x_{ji} \leq 1, x_{ji} \geq 0.$$

In this formulation, the number of constraints is increased by $\sum_i |\mathcal{P}_i|$ which can be kept constant independent of the number of impressions and users.

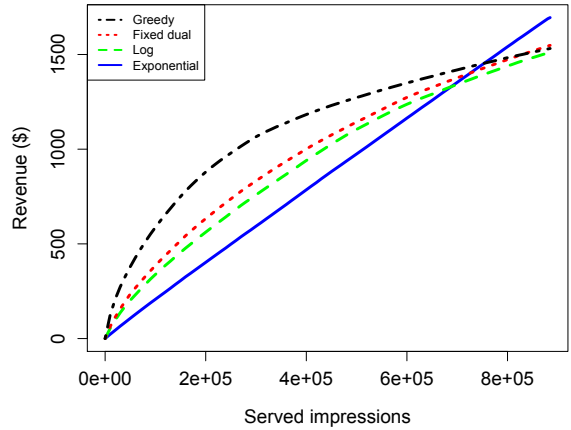

**Figure 8** Revenue under fcap constraints.

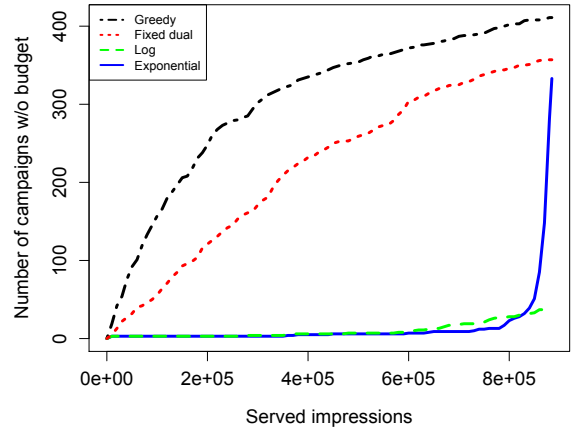

**Figure 9** Number of campaigns with no budget under fcap constraints.

**Table 3** Performance of allocation algorithms with Fcap.

| Allocation algorithm | Total Revenue | Improvement over greedy | Mid flight oob | Final oob |
|---|---|---|---|---|
| Greedy | $1531.71 | 0.0% | 335 | 411 |
| Fixed dual | $1548.32 | 1.08% | 232 | 357 |
| Log | $1511.20 | -1.3% | 6 | 39 |
| Exponential | $1694.89 | 10.6% | 5 | 333 |

Similar to Figures 6 and 7, revenue and budget delivery of the four allocation algorithms are represented in Figures 8 and 9. For this experiment, we have used 10 quintiles of non-zero bids of each advertiser as a partition on bids. Again, exponential value allocation gives the best revenue and most balanced delivery. Because of using partitions instead of actual frequency capping constraints, estimated dual prices are less accurate. Thus, it is not a surprise to observe a smaller improvement over greedy in the dual price allocation algorithm.

# 7. Conclusion, Discussion, and Future Work

In this paper, we developed dual price updating rules that incorporate the belief of the ad network toward the future demand. We modeled the impression allocation problem as a convex risk minimization problem, and impressions are allocated such that the exposed risk is minimized. Such a convex representation leads us to build a connection between the belief of the ad network toward the future demand and dual price updating rules. Price updating rules have simple representations in terms of the remaining budgets, and they can be implemented in a real ad network with billions of impressions served per day. The proposed allocation algorithms are evaluated for a real data set from AOL/Advertising.com. We empirically showed that the exponential updating rule has the best performance in terms of the total collected revenue and the smooth budget consumption.

Even though the empirical results are presented by applying the algorithms in two consequent time periods ($T_1$ and $T_2$), they can be extended to a rolling horizon scheme. In Section 6, initial dual prices are estimated using a sample of impressions from the first time period ($T_1$). The estimated dual prices are used in the second time period $T_2$ to evaluate the performance. However, when there are $d$ time periods, $T_1, T_2 \ldots T_d$, the allocation algorithm in each time period utilizes the dual prices that are obtained from the pervious time period.

One direction for future work is addressing the delay in the observation of revenue and cost information. In practice, a time delay in observing the total exhausted budget or the collected revenue is inevitable. This can be due to the nature of the revenue model, or it can be due to technological limitations. We observed that the pro- posed algorithms are very sensitive to the delay, and their performance degrades considerably by introducing even a small amount of delay.

In summary, we proposed simple price updating rules that minimized the ex- posed risk to the ad network. These updating rules perform well on real data sets without assuming any structure on the sequence of impression arrivals.

**Appendix A: Monotonicity of** $\theta(.)$

LEMMA 6. *If* $\forall s > 0, \frac{\partial v_i}{\partial s}(s) \leq p_{max}$ *then*

$$X_S \geq X_{\bar{S}} \ \Rightarrow \ \theta(X_S) \leq \theta(X_{\bar{S}}).$$

*Proof* For all $\mathbf{p} \in \Omega$, $X_S(\mathbf{p}) \geq X_{\bar{S}}(\mathbf{p})$, or equivalently,

$$w + \sum_{i=1}^{N} p_i s_i \geq \bar{w} + \sum_{i=1}^{N} p_i \bar{s}_i, \tag{36}$$

$$w - \bar{w} \geq \sum_{i=1}^{N} p_i (\bar{s}_i - s_i). \tag{37}$$

On the other hand,

$$\theta(X_{\bar{S}}) - \theta(X_S) = -\bar{w} + w + \sum_{i=1}^{N} v_i(s_i) - v_i(\bar{s}_i) \tag{38}$$

$$= -\bar{w} + w + \sum_{i=1}^{N} \frac{\partial v_i}{\partial s}(\hat{s}_i)(s_i - \bar{s}_i),$$

where $\hat{s}_i$ is a real value between $s_i$ and $\bar{s}_i$. Since equation (37) is satisfied for all $0 \leq p_i \leq p_{\max}$, from equations (37) and (38), we can conclude that $\theta(X_{\bar{S}}) - \theta(X_S) \geq 0$. Q.E.D.

**Appendix B: Proof of Lemma 4**

Using Theorem 4.30 in Föllmera and Schied (2004), we just need to show $\rho(.)$ is continuous from above, i.e.,

$$X_n \searrow X \Longrightarrow \rho(X_n) \nearrow \rho(X).$$

Assume $X_n \searrow X$. By the monotonicity property of the risk measure $\rho(.)$,

$$\rho(X) \geq \rho(X_{n+1}) \geq \rho(X_n).$$

Thus, $\rho(X) \geq \lim_{n \to +\infty} \rho(X_n)$.

For every $\epsilon \geq 0$, and for every $n$, there exists $X_{S_n} \in \hat{\mathcal{X}}$ such that

$$X_{S_n} \leq X_n, \tag{39}$$

$$\theta(X_{S_n}) \leq \rho(X_n) + \epsilon. \tag{40}$$

Consider sequence $\{S_n\}_i$. Since $0 \leq X_{S_n} \leq X_n \leq X_1$ and $X_1$ is bounded by definition, this sequence is bounded in $\mathbb{R}^{M+1}$. Thus there is a subsequence $\{S_{n'}\}$ that converge to a point $S^*$. Since $X_{S_{n'}} \leq X_{n'}$,

$$X_{S^*} = \lim_{n' \to +\infty} X_{S_{n'}} \leq \lim_{n' \to +\infty} X_{n'} = X.$$

Since $v_i(.)$ are continuous, thereby,

$$\rho(X) \leq \theta(X_{S^*}) \leq \lim_{n \to +\infty} \rho(X_n) + \epsilon. \tag{41}$$

Since equation (41) is valid for every $\epsilon > 0$, $\rho(X) \leq \lim_{n \to +\infty} \rho(X_n)$, and consequently,

$$\rho(X) = \lim_{n \to +\infty} \rho(X_n).$$

**Appendix C: Proof of Lemma 5**

Recall, by Lemma 4, $\alpha(Q)$ can be restated by

$$\alpha(Q) = \sup_{X \in \mathcal{X}} E_Q[-X] - \rho(X). \tag{42}$$

First, we show that, in the above representation, taking supremum over $\hat{\mathcal{X}}$ is equivalent of taking supremum over $\mathcal{X}$. Fix a positive real value $\epsilon > 0$. By the definition of risk measure (5), for every $X \in \mathcal{X}$, there exists $X_S^\epsilon (\in \hat{\mathcal{X}}) \leq X$ such that

$$\rho(X) \leq \theta(X_S^\epsilon) = \rho(X_S^\epsilon) \leq \rho(X) + \epsilon. \tag{43}$$

Therefore,

$$E_Q[-X] - \rho(X) \leq E_Q[-X_S^\epsilon] - \rho(X_S^\epsilon) + \epsilon.$$

Hence,

$$\begin{aligned}\alpha(Q) &= \sup_{X \in \mathcal{X}} E_Q[-X] - \rho(X) \\ &\leq \sup_{X_S \in \hat{\mathcal{X}}} E_Q[-X_S] - \rho(X_S) + \epsilon.\end{aligned} \tag{44}$$

Since equation (44) is valid for every $\epsilon > 0$, thus

$$\alpha(Q) \leq \sup_{X_S \in \hat{\mathcal{X}}} E_Q[-X_S] - \rho(X_S). \tag{45}$$

On the other hand, since $\hat{\mathcal{X}} \subset \mathcal{X}$

$$\sup_{X_S \in \hat{\mathcal{X}}} E_Q[-X_S] - \rho(X_S) \leq \alpha(Q). \tag{46}$$

Therefore, from (45) and (46)

$$\alpha(Q) = \sup_{X_S \in \hat{\mathcal{X}}} E_Q[-X_S] - \rho(X_S). \tag{47}$$

Next, observe

$$E_Q[-X_S] = E_Q[-w - \sum_{i=1}^N p_i s_i] \tag{48}$$

$$= -w - \sum_{i=1}^N E_Q[p_i] s_i \tag{49}$$

$$= -w - \sum_{i=1}^N p_i^* s_i. \tag{50}$$

And finally, since $\rho(X_S) = -w - \sum_{i=1}^{N} v_i(s_i)$,

$$\begin{aligned}
\alpha(Q) &= \sup_{X_S \in \hat{\mathcal{X}}} E_Q[-X_S] - \rho(X_S) \\
&= \sup_{X_S \in \hat{\mathcal{X}}} -\sum_{i=1}^{N} p_i^* s_i + \sum_{i=1}^{N} v_i(s_i) \\
&= \sup_{X_S \in \hat{\mathcal{X}}} \sum_{i=1}^{N} \Big(v_i(s_i) - p_i^* s_i\Big).
\end{aligned} \quad (51)$$

Q.E.D.

## Appendix D: Logarithmic Risk Measure

In this appendix, we compute probability penalty function $\alpha(Q)$ for risk measure in (6) when

$$v_i(s_i) = p_i^\epsilon s_i + \gamma \log\left(\frac{s_i}{b_i}\right).$$

Substituting $v_i(s_i)$ in (17), we have

$$\alpha(Q) = \sup_{X \in \mathcal{X}} \sum_{i=1}^{N} \left( p_i^\epsilon s_i + \gamma \log\left(\frac{s_i}{b_i}\right) - p_i^* s_i \right). \quad (52)$$

Note that since $s_i$s are independent, each term of equation (52) can be maximized separately. Let $\beta_i(s_i) = p_i^\epsilon s_i + \gamma \log\left(\frac{s_i}{b_i}\right) - p_i^* s_i$. Then,

$$\frac{\partial \beta_i}{\partial s_i} = p_i^\epsilon + \frac{\gamma}{s_i} - p_i^*, \quad (53)$$

and

$$\frac{\partial^2 \beta_i}{\partial s_i^2} = -\frac{\gamma}{s_i^2} < 0. \quad (54)$$

Since the second derivative is negative, $\beta_i(s_i)$ achieves its maximum when $\beta_i(s_i) = 0$, or on the boundaries. There are two cases.

(a) If $b_i(p_i^* - p_i^\epsilon) \leq \gamma$ then $\frac{\partial \beta_i}{\partial s_i} \leq 0$ and $\beta_i(s_i)$ achieves its maximum on the boundary $s_i^* = b_i$.

(b) If $b_i(p_i^* - p_i^\epsilon) > \gamma$ function $\beta_i(s_i)$ is maximized at point $s_i^*$ where $\frac{\partial \beta_i}{\partial s_i}(s_i^*) = 0$, or equivalently,

$$p_i^* = p_i^\epsilon + \frac{\gamma}{s_i^*}, \quad (55)$$

$$s_i^* = \frac{\gamma}{p_i^* - p_i^\epsilon}. \quad (56)$$

Combining (a) and (b), $\beta_i(s_i)$ is maximized at

$$s_i^* = \begin{cases} b_i & b_i(p_i^* - p_i^\epsilon) \leq \gamma \\ \frac{\gamma}{p_i^* - p_i^\epsilon} & \text{o.w.} \end{cases}$$

If $b_i(p_i^* - p_i^\epsilon) \leq \gamma$,
$$\max_{s_i} \beta_i(s_i) = b_i p_i^\epsilon - b_i p_i^*.$$

If $b_i(p_i^* - p_i^\epsilon) > \gamma$,
$$\max_{s_i} \beta_i(s_i) = p_i^\epsilon s_i^* + \gamma \log\left(\frac{s_i^*}{b_i}\right) - p_i^* s_i^* \tag{57}$$

$$= p_i^\epsilon \frac{\gamma}{p_i^* - p_i^\epsilon} + \gamma \log\left(\frac{\frac{\gamma}{p_i^* - p_i^\epsilon}}{b_i}\right) - p_i^* \frac{\gamma}{p_i^* - p_i^\epsilon} \tag{58}$$

$$= -\gamma + \gamma \log\left(\frac{\gamma}{b_i p_i^* - b_i p_i^\epsilon}\right). \tag{59}$$

Or,
$$\alpha(Q) = \sum \beta_i(s_i) = \|\mathbf{p}^\epsilon.\mathbf{b} - \mathbf{p}^*.\mathbf{b}\|_\lambda,$$

where function $\|\mathbf{x}\|_\lambda$ is defined in (22).

### Appendix E: Exponential Risk Measure

Here we derive the probability penalty function for the exponential risk measure. In the exponential risk measure, the budget value function, $v_i(.)$, is defined as follows.

$$v_i(s_i) = p_i^\epsilon \left(-\frac{b_i}{\kappa}\right) e^{k\left(\frac{b_i - s_i}{b_i} - \frac{h}{N}\right)}.$$

By Lemma 5,
$$\alpha(Q) = \sup_{X \in \mathcal{X}} \sum_{i=1}^{N} \left( p_i^\epsilon \left(-\frac{b_i}{\kappa}\right) e^{k\left(\frac{b_i - s_i}{b_i} - \frac{h}{N}\right)} - p_i^* s_i \right). \tag{60}$$

Similar to the previous risk measures, each term in (60) can be maximized independently. Let
$$\alpha_i(s_i) = p_i^\epsilon \left(-\frac{b_i}{\kappa}\right) e^{k\left(\frac{b_i - s_i}{b_i} - \frac{h}{N}\right)} - p_i^* s_i. \tag{61}$$

In order to find the optimal solution, we the find first and the second derivative of $\alpha_i(s_i)$.

$$\frac{\partial \alpha_i}{\partial s_i} = p_i^\epsilon e^{k\left(\frac{b_i - s_i}{b_i} - \frac{h}{N}\right)} - p_i^*. \tag{62}$$

$$\frac{\partial^2 \alpha_i}{\partial s_i^2} = p_i^\epsilon \left(-\frac{\kappa}{b_i}\right) e^{k\left(\frac{b_i - s_i}{b_i} - \frac{h}{N}\right)} < 0. \tag{63}$$

Since the second derivative is negative, $\alpha_i(s_i)$ attains its maximum where $\frac{\partial \alpha_i}{\partial s_i}(s_i^*) = 0$ or on the boundaries. Thus,

$$p_i^* = p_i^\epsilon e^{\kappa\left(\frac{b_i - s_i^*}{b_i} - \frac{h}{N}\right)} \tag{64}$$

$$\frac{p_i^*}{p_i^\epsilon} = e^{\kappa\left(\frac{b_i - s_i^*}{b_i} - \frac{h}{N}\right)} \tag{65}$$

$$s_i^* = b_i \left(1 - \frac{1}{\kappa} \log\left(\frac{p_i^*}{p_i^\epsilon}\right) - \frac{h}{N}\right). \tag{66}$$

Substituting $s_i^*$ in equation (61),

$$\alpha_i(s_i^*) = p_i^\epsilon \left(-\frac{b_i}{\kappa}\right) e^{k\left(\frac{b_i - s_i^*}{b_i} - \frac{h}{N}\right)} - p_i^* s_i^* \tag{67}$$

$$= p_i^\epsilon \left(-\frac{b_i}{\kappa}\right) \frac{p_i^*}{p_i^\epsilon} - p_i^* b_i \left(1 - \frac{1}{\kappa} \log\left(\frac{p_i^*}{p_i^\epsilon}\right) - \frac{h}{N}\right) \tag{68}$$

$$= p_i^* b_i \left(-\frac{1}{\kappa} - 1 + \frac{1}{\kappa} \log\left(\frac{p_i^*}{p_i^\epsilon}\right) + \frac{h}{N}\right). \tag{69}$$

Finally, the probability penalty function $\alpha(Q)$ can be represented by KL-divergence as follows.

$$\alpha(Q) = \sum_{i=1}^N \alpha_i(s_i^*) \tag{70}$$

$$= \sum_{i=1}^N p_i^* b_i \left(-\frac{1}{\kappa} - 1 + \frac{1}{\kappa} \log\left(\frac{p_i^*}{p_i^\epsilon}\right) + \frac{h}{N}\right) \tag{71}$$

$$= \left(-\frac{1}{\kappa} - 1 + \frac{h}{N}\right) \sum_{i=1}^N p_i^* b_i + \frac{1}{\kappa} \sum_{i=1}^N p_i^* b_i \left(\log\left(\frac{p_i^*}{p_i^\epsilon}\right)\right) \tag{72}$$

$$= \sum_{i=1}^N p_i^* b_i \left(-\frac{1}{\kappa} - 1 + \frac{h}{N} + \frac{1}{\kappa} \sum_{i=1}^N \frac{p_i^* b_i}{\sum_{j=1}^N p_j^* b_j} \log\left(\frac{p_i^* b_i}{p_i^\epsilon b_i}\right)\right) \tag{73}$$

$$= \left(\sum_{i=1}^N p_i^* b_i\right) \left(-\left(1 - \frac{h}{N}\right) + \frac{1}{\kappa}\left(\mathcal{D}_{KL}\left(\mathbf{p}^\epsilon.\mathbf{b} \| \mathbf{p}^*.\mathbf{b}\right) - 1\right)\right). \tag{74}$$

## Acknowledgments

The authors gratefully acknowledge Bin Ren who generously offered help in reviewing the theoretical parts of the paper. The fourth author is partially supported by AFOSR Grant FA9550-12-1-0396.